\documentstyle[12pt]{amsart}

\begin{document}
\vspace*{.7cm}
\begin{center}{\LARGE \bf MODULAR FORMS ASSOCIATED WITH THE MONSTER
 MODULE}
\end{center}
\vspace{.4cm}
\begin{center}
\begin{tabular}{cc}
\mbox{Koichiro Harada}  & Mong Lung Lang\\
Department of Mathematics & Department of Mathematics\\
The Ohio State University  & National University of Singapore\\
Columbus, Ohio 43210 USA & Republic of Singapore\\
\end{tabular}
\end{center}

\vspace{.6cm}

\baselineskip=24pt


\section{Introduction}

In Harada-Lang [4], we associated to each irreducible character $\chi$
 of the monster simple group $ \Bbb M$ a modular function $t_{\chi}(z)$,
 called in [4], the McKay-Thompson series for $\chi$.
 $t_{\chi}(z)$ is a weighted average of all McKay-Thompson series
 $t_g(z)$ for the element $g$ of $\Bbb M$ as $g$ ranges over $\Bbb M$ :
$$ t_{\chi}(z) = \frac{1}{|\Bbb M|} \sum_{g \in\Bbb M} \chi(g) t_g(z).$$
If $\Gamma_{\chi}$ is the invariance subgroup of $t_{\chi}(z)$, then we
 showed
$$\Gamma_{\chi} = \Gamma_0(N_{\chi}) = \bigcap_{g \in\Bbb M} \Gamma_g$$
where $g$ ranges over all the elements of $\Bbb M$ such that
$\chi(g) \ne 0$ and
 $$ N_{\chi}  =  \text {lcm}\{ n_gh_g : \, \text{for all} \, g \in \ \Bbb M
 \, \text{with}  \,  \chi(g) \ne 0 \}.$$
As shown in Conway-Norton [1], the invariance group $\Gamma_g$
 of $t_g(z)$ is a certain subgroup of index $h$ of the conjugate by
 $$\left(\begin{array}{ll}
                h & 0 \\
                  0  & 1 \\
           \end{array}
    \right) $$
of
$$ \Gamma_0(\frac{n}{h}) + e, f, \cdot  \cdot  \cdot  $$
where $e$, $f,$ etc. denote the Atkin-Lehner involutions.
 In [1], such a conjugate is denoted by
$$ n|h + e, f,\cdot \cdot \cdot \, .$$
The numbers
 $n$, $h$ depend on $g$, hence our notation $n_g$, $h_g$.
  Obviously every  $t_{\chi}(z)$ is  invariant by
$$\bigcap_{g \in \Bbb M} \Gamma_g = \Gamma_0(N_0)$$
 where $N_0$ = $2^63^35^27\cdot 11 \cdot 13 \cdot 17 \cdot
 19 \cdot 23 \cdot 29 \cdot 31\cdot 41 \cdot 47 \cdot 59 \cdot 71
\sim 10^{21}.$
The level $N_{\chi}$ can be very large or relatively small. For example,
$$N_{\chi_1} = N_0, \, N_{\chi_{166}} = 2^63^37 = 4032$$
where $\chi_1$ = $1$ is the trivial character and the character
 numbering such as $\chi_{166}$ is taken from the Atlas.
  In this paper, we will investigate the relation between
 $t_{\chi}(z)$ and the generating functions of the highest weight vectors
 (also called singular vectors, primary fields or lowest weight vectors.)

\bigskip

\section{the monster module as a $Vir$ module}

The monster module $\Bbb V$ is constructed in Frenkel-Lepowsky-Meurman
 [3] as a vertex operator algebra and is denoted by
  $\Bbb V^{\natural}$
 there.
  Let $V$ be a vertex operator algebra.
  Then $V$ possesses two distinguished elements $1$ and $\omega$, called
the vacuum and the conformal vector (or the Virasoro element)
 of $V$, respectively.

If $Y(\omega, z) = \sum \omega_nz^{-n-1}$ is the vertex operator
 corresponding to the conformal vector $\omega$ and if
 we set $L(n) = \omega_{n+1}$ for $ n \in \Bbb Z$, then $L(n)$
 satisfies the commutation relation:
$$[L(n), L(m)] = (n-m)L_{n+m} + \frac{1}{12}(n^3 - n)c\delta_{n+m,0}$$
where $c$ is a constant called the central charge of $V$.
For the monster module $\Bbb V$, $c$ = 24. $c$ is also called the rank
 of the vertex operator algebra $V$.

Let $\cal L$ be the Lie algebra generated by $L(n)$, $n \in \Bbb Z.$
 $\cal L$ is denoted by $Vir$ else where.  The subalgebras
 $\cal L^+$ and $\cal L^-$ are generated by $L(n), n \in \Bbb Z^+$
 and  $L(n), n \in \Bbb Z^-$, respectively.  It is known
 that $\Bbb V$ possesses a positive definite invariant
 bilinear form and so $\Bbb V$ is completely  reducible as an $\cal L$
module and is a sum of highest weight modules.

Let $M(h,c)$ be the Verma  module of the Virasoro algebra of
 central charge $c$ generated by the highest weight vector $v$ of
 height $h$ : i.e.
$$ M(h,c) = \cal L v, \, \cal L^+ v = 0, \, L(0)v = hv.$$

The module structure of $M(h,c)$ has been determined by
 Feigin-Fuchs [2].  We will use their results to determine the module
 structure of $\Bbb V$ as an $\cal L$ module.  Feigin-Fuchs showed
 that every submodule of $M(h,c)$ is a sum of submodules that are
 also Verma modules.  Therefore, the knowledge of all embeddings
 among Verma modules gives all submodules of a given Verma module.
  The main theorem of Feigin-Fuchs states that there are six types
 of embeddings of the Verma modules into
 other Verma modules. Let
$$\cases  p\alpha - q\beta = m \\
        c = 24 = \frac{(3p - 2q)(3q-2p)}{pq} \\
       h = \frac{m^2 - (p-q)^2}{4pq}  \endcases   \eqno (1) $$
where $p$, $q$ and $m$ are complex variables. We next solve for
integers $\alpha$ and $\beta$. Let
$$\epsilon = \frac{-11 \pm i \sqrt {23}}{2}, \quad
  \bar \epsilon = \frac{-11 \mp i \sqrt {23}}{2}$$
We compute
$$ \epsilon \bar\epsilon = 1, \,
  \epsilon +  \bar\epsilon = \frac{-11}{6},
 \,  \epsilon^2 +  \bar\epsilon^2 = \frac{49}{36}.$$
Using the second equality of (1), we obtain
$$(p\alpha -q\beta)^2 = m^2 = 4pq + (q-p)^2,$$
which may be rewritten as
$$(\alpha - \epsilon \beta)^2 = 4 \epsilon h + (\epsilon -1)^2.$$
We therefore obtain two equations :
$$\alpha^2 - 2\epsilon \alpha \beta + \epsilon^2 \beta^2 =
  4\epsilon h + (\epsilon - 1)^2,$$
and
$$\alpha^2 - 2\bar \epsilon \alpha \beta + \bar \epsilon^2 \beta^2 =
  4\bar \epsilon h + (\bar \epsilon - 1)^2.$$
Taking the sum of them, we get
$$ 72\alpha^2 + 132\alpha\beta +49\beta^2 = -264h + 253.$$
By subtracting one from the other, we get
$$-12\alpha\beta - 11 \beta^2 = 24h - 23.$$
 Therefore
$$\alpha^2 - \beta^2 = 0,$$
or $\alpha = \pm \beta.$  Setting $\alpha = \delta\beta$ with
 $\delta = \pm 1$, we have
$$ \beta^2 = \frac{24h-1}{11-12\delta}.$$
If $h$ = $0$, then we must have $\delta = 1$ and so
$ \alpha = \beta = \pm 1.$  In particular, $ \alpha \beta = 1 > 0$.
  On the other hand, if $h \in \Bbb Z^+$, then
 $\delta = -1$ and so
$ \alpha = - \beta = \pm 1, $ and hence $\alpha \beta = -1 < 0 .$
  Using the results of Feigin-Fuchs [2], we conclude
(which must be well known to experts) :

{\bf Theorem. } {\em $M(0,24)$ has a unique submodule, which
 is isomorphic to  $M(1,24)$.  For all positive integers $h$,
 $M(h,24)$ is irreducible. }

Let $L(c,h)$ be the unique irreducible highest weight $\cal L$-module
 of central charge $c$ and height $h$.  Then

{\bf Corollary.}  {\em We have}

$(1)$. $L(0,24) = M(0,24)/M(1,24),$ {\em and,}

$(2).$ $L(h,24) = M(h, 24)$ {\em if} $h \in \Bbb Z^+.$

Let us now express the monster module $\Bbb V$ as a sum
 of $ L(h,24)$'s as follows
$$\Bbb V = \sum_{h = 0}^{\infty}s_h L(h,24).$$
Then $s_h$ is the number of linearly independent singular
 vectors $v_h$ of height $h$, hence $v_h \in \Bbb V_h$.  Since
 the Virasoro algebra $\cal L$ commutes with the action of the monster
 $\Bbb M$, we can actually split $s_h$ into the sum of $s_h^k$ where the
index $k$ corresponds to the irreducible character $\chi_k$. More precisely,
 let
$$\Bbb V_h^k = c_{hk}\chi_k$$
where $ c_{hk}$ is the multiplicity of $\chi_k$
 in $\Bbb V_h$ and
$$\Bbb V^k = \coprod_{h = 0}^{\infty} \Bbb V_h^k.$$
Thus $\Bbb V^k$ is an $\Bbb M$ submodule of $\Bbb V$ consisting
 entirely of irreducible submodules isomorphic to
 $\chi_k$ and $\Bbb V_h^k$ is an $\Bbb M$ submodule of $\Bbb V^k$
 of height $h$.  We also define
$$W_h^k = \cal L (\coprod_{0 \le i < h} \Bbb V_i^k) \cap \Bbb V_h^k,$$
which is an $\Bbb M$ submodule of $\Bbb V_h^k$ that is generated by
 elements of lower heights.  Let
$$ s_h^k = \text{dim}\Bbb V_h^k/W_h^k.$$
Then $s_h^k$ is the number of linearly independent singular vectors in
$\Bbb V_h^k$. Obviously
$$ s_h = \sum_{k = 1}^{194} s_h^k.$$
For a graded module $X = \sum_{h \in \Bbb Z} X_h$, the
 character of $X$ (or the partition function of $X$)
 is defined to be  a formal sum
$$\text{char}(X) = \sum_{h \in \Bbb Z} \text{dim} X_h x^h.$$
Using this notation, we have, as is well known,
$$\text{char}M(h,c) = x^h \sum_{n \ge 0} p(n)x^n$$
where $p(n)$ is the partition function of $n$.  For convenience,
 set $p(0) = 1$, and $p(n) = 0$ if $n \in \Bbb Z^-.$  Let us
 consider the $\cal L$ submodule generated
 by the vacuum $1$.  We have $V_1 = 0$ but the height 1 component of
$M(0,24)$ is nonzero, we conclude that
$$\cal L\cdot 1 \simeq M(0,24)/M(1,24).$$
Hence
$$\text{char}(\cal L\cdot 1) = \sum_{n \ge 0} p(n)x^n -
 x \sum_{n \ge 0} p(n)x^n
               = \sum_{n \ge 0}(p(n) - p(n-1))x^n.$$
Writing
$$\text{char}(\cal L\cdot 1) = \sum_{h \ge 0} a_{h1} x^h,$$
 we get a partial list :
$$\begin{array}{llllllllllll}
h & 0 & 2 &3 & 4 &5 & 6 & 7 & 8 & 9 & 10 &11 \\
a_{h1}&1 &1 & 1 & 2& 2 & 4 &4 &  7& 8 & 12 & 14 \\
\end{array}  $$
In [4], we had a partial list of $c_{h1}$ where
  $c_{h1}$ is the multiplicity of the trivial
 character $\chi_1$ occuring in $\Bbb V_h$.
$$\begin{array}{llllllllllll}
h & 0 & 2 &3 & 4 &5 & 6 & 7 & 8 & 9 & 10 &11  \\
c_{h1}&1 &1 & 1 & 2& 2 & 4 &4 &  7& 8 & 12 & 14 \\
\end{array}  $$
The  coincidence  $ c_{h1} = a_{h1}$ stops there and we have
$$\begin{array}{ll}
      h &  12 \\
     a_{h1} & 21 \\
     c_{h1} & 22  \\
   \end{array}$$
This means $s_{12}^1 = 1$, namely, $\Bbb V_{12}^1$ contains a singular
 vector, while $\Bbb V_h^1$, $0 < h \le 11$, do not.
 The number $d$ of linearly independent singular vectors occuring
  in $\Bbb V_h^1$ ($0 \le h \le 30$) is as follows
$$\begin{array}{rrrrrrrrrrrr}
h & 12 & 16 &18 & 20 &22 & 24 & 26 & 27 & 28 & 29&30  \\
d &1 &1 & 1 & 1& 1&3 &2 &  1& 4 & 2 & 6  \\
\end{array}  $$
We are now lead to consider its generating function for
each $k$, $1 \le k \le 194$. Define
$$G^k(x) = \sum_{h \ge 0} s_h^k x^h.$$
The character of $\Bbb V^k$ is
$$\text{char}(\Bbb V^k) = \sum_{h \ge 0} c_h^k(\text{deg}\chi_k)x^h
 = x \text{deg}\chi_k t_{\chi}(x)$$
 where $t_{\chi}(z)$ is the McKay-Thompson series for the
 irreducible character $\chi$.  On the other hand, using the expression
$$\Bbb V^k = \sum_{h \ge 0} s_h^kL(h,24),$$
we obtain
$$\text{char}(\Bbb V^k) =  \sum_{h \ge 0} s_h^k \text{char}L(h,24).$$
Suppose $k > 0$. Then $s_0^k = 0$ and so
$$ \text{char}(\Bbb V^k) =
 \sum_{h \ge 1} s_h^k x^h \sum_{n \ge 0} p(n)x^n.$$
 On the other hand if $k = 1$, then $L(0,24)$
 occurs only once as a constitient of $\Bbb V^1$. Therefore
$$\text{char}(\Bbb V^1) = (1 - x + \sum_{h \ge 2}s_h^1x^h)\sum_{n \ge 0}
p(n)x^n.$$
Using the Dedekind eta-function and replacing $x$ by
 $q = e^{2 \pi i z}$, we obtain, by setting
 $s_1^1 = -1$  for convenience,
$$\text{deg} \chi_k t_{\chi_k}(q) =
\frac{q^{-1}(\sum_{h \ge 0} s_h^k q^h)q^{\frac{1}{24}}}{\eta ( q)}.$$
Hence
$$\text{deg} \chi_k t_{\chi_k}(q) \eta(q) = q^{-\frac{23}{24}} \sum_{h \ge 0}
s_h^k q^h,$$
which implies
$$q^{-\frac{23}{24}}G^k(q) = \text{deg}\chi_k t_{\chi_k}(q)\eta(q)$$
where as defined before $G^k(q)$ is the generating function
 of the singular vectors in $\Bbb V^k$.
 Writing $G^k = G^{\chi}$ in general, we obtain :

{\bf Theorem.}  $q^{-\frac{23}{24}} G^{\chi}(q)$ {\em is a meromorphic
 modular form of weigh}t $\frac{1}{2}$ {\em and level} $N_{\chi}.$

{\bf Corollary.} $ q^{-\frac{23}{24}} G^{\chi}(q)
\eta(q)^{23}$ {\em is a holomorphic modular function of weight $12$ and
 level $N_{\chi}.$ }

The first 52 coefficients of $\frac{G^{\chi _i}(q)}{{\text deg}\chi_i},$ $ 1
\le i \le 194,$
 can be found in the following
 table.

\newpage
\small
$
\begin{array}{r}
a_n \\a_0\\a_1\\a_2\\a_3\\a_4\\a_5\\a_6\\a_7\\a_8\\a_9\\
a_{10}\\a_{11}\\a_{12}\\a_{13}\\a_{14}\\a_{15}\\a_{16}\\a_{17}\\a_{18}\\a_{19}\\
a_{20}\\a_{21}\\a_{22}\\a_{23}\\a_{24}\\a_{25}\\a_{26}\\a_{27}\\a_{28}\\a_{29}\\
a_{30}\\a_{31}\\a_{32}\\a_{33}\\a_{34}\\a_{35}\\a_{36}\\a_{37}\\a_{38}\\a_{39}\\
a_{40}\\a_{41}\\a_{42}\\a_{43}\\a_{44}\\a_{45}\\a_{46}\\a_{47}\\a_{48}\\a_{49}\\
a_{50}\\a_{51}\\
\end{array}$
$\begin{array}{r}
\chi_{1}\\
1\\-1\\0\\0\\0\\0\\0\\0\\0\\0\\0\\0\\1\\0\\0\\0\\1\\0\\1
\\0\\1\\0\\1\\0\\3\\0\\2\\1\\4\\2\\6\\3\\8\\5\\10\\8\\17\\11\\
21\\19\\31\\27\\45\\40\\63\\62\\89\\91\\133\\132\\185\\200\\
\end{array}$
$\begin{array}{r}
\chi_{2}\\
0\\0\\1\\0\\0\\0\\0\\0\\0\\0\\1\\0\\1\\1\\2\\1\\3\\2\\4\\
3\\6\\5\\9\\8\\13\\14\\21\\21\\31\\34\\48\\53\\73\\82\\112\\
128\\170\\198\\261\\304\\396\\470\\604\\720\\920\\1102\\1402\\
1688\\2129\\2579\\3244\\3932\\
\end{array}$
$\begin{array}{r}
\chi_{3}\\
0\\0\\0\\1\\0\\0\\0\\0\\0\\1\\0\\1\\1\\2\\1\\3\\2\\4\\4\\6\\
6\\11\\10\\16\\18\\25\\28\\41\\45\\63\\74\\99\\116\\158\\184\\
245\\296\\383\\463\\604\\729\\940\\1151\\1468\\1800\\2296\\2817
\\3573\\4412\\5566\\6878\\8677\\
\end{array}$
$\begin{array}{r}
\chi_{4}\\
0\\0\\0\\0\\1\\0\\0\\0\\1\\0\\1\\1\\2\\1\\3\\2\\5\\4\\7\\8\\
13\\12\\20\\22\\33\\37\\53\\60\\87\\100\\137\\165\\224\\266\\
358\\435\\574\\702\\919\\1130\\1473\\1817\\2344\\2914\\3746\\
4653\\5962\\7436\\9483\\11851\\15076\\18858\\
\end{array}$
$\begin{array}{r}
\chi_{5}\\
0\\0\\0\\0\\0\\0\\1\\0\\1\\0\\1\\0\\2\\1\\3\\2\\5\\4\\9\\7\\
14\\14\\23\\24\\39\\41\\64\\72\\105\\121\\175\\204\\286\\343\\
469\\569\\771\\942\\1259\\1556\\2054\\2553\\3352\\4181\\5451\\
6835\\8859\\11137\\14390\\18120\\23330\\29450\\
\end{array}$
$\begin{array}{r}
\chi_{6}\\
0\\0\\0\\0\\0\\1\\0\\1\\0\\1\\1\\2\\1\\3\\3\\5\\5\\9\\9\\15\\
16\\25\\29\\42\\48\\70\\83\\115\\139\\192\\233\\315\\388\\517\\
643\\849\\1057\\1388\\1739\\2262\\2845\\3687\\4642\\5991\\7558\\
9719\\12283\\15752\\19914\\25493\\32262\\41208\\
\end{array}$
$\begin{array}{r}
\chi_{7}\\
0\\0\\0\\0\\0\\0\\1\\0\\0\\1\\1\\1\\2\\2\\3\\4\\5\\7\\11\\12\\
17\\23\\31\\39\\54\\68\\92\\119\\156\\202\\269\\342\\450\\583\\
761\\978\\1278\\1641\\2134\\2747\\3555\\4574\\5922\\7599\\9811\\
12611\\16249\\20859\\26864\\34466\\44333\\56884\\
\end{array}$
$\begin{array}{r}
\chi_{8}\\
0\\0\\0\\0\\0\\0\\0\\1\\1\\1\\1\\2\\2\\4\\4\\7\\9\\13\\15\\24\\
30\\42\\53\\75\\95\\132\\169\\230\\297\\400\\514\\692\\893\\1185\\
1534\\2029\\2623\\3456\\4470\\5861\\7588\\9918\\12826\\16737\\21644
\\28162\\36418\\47318\\61155\\79358\\102544\\132903\\
\end{array}$
$\begin{array}{r}
\chi_{9}\\
0\\0\\0\\0\\0\\0\\0\\0\\1\\0\\1\\1\\2\\2\\4\\4\\8\\8\\14\\16\\26\\31\\48\\59\\87\\108\\156\\197\\279\\356\\494\\636\\872\\1126\\1528\\1982\\2667\\3467\\4637\\6034\\8029\\10460\\13857\\18065\\23855\\31102\\40962\\53418\\70199\\91549\\120111\\156644\\
\end{array}$
$\begin{array}{r}
\chi_{10}\\
0\\0\\0\\0\\0\\0\\0\\0\\0\\1\\1\\1\\2\\2\\4\\5\\7\\10\\15\\20\\28\\38\\53\\71\\99\\133\\182\\244\\331\\443\\600\\797\\1072\\1427\\1910\\2532\\3379\\4472\\5947\\7862\\10421\\13755\\18201\\23984\\31671\\41701\\54973\\72303\\95204\\125119\\164575\\216162\\
\end{array}$
\newpage
$\begin{array}{r}
\chi_{11}\\
0\\0\\0\\0\\0\\0\\0\\0\\0\\0\\1\\0\\1\\0\\2\\1\\3\\3\\7\\6\\12\\13\\23\\25\\42\\49\\77\\93\\140\\174\\255\\320\\457\\586\\821\\1058\\1464\\1899\\2599\\3389\\4597\\6013\\8109\\10623\\14248\\18707\\24988\\32832\\43727\\57497\\76390\\100508\\
\end{array}$
$\begin{array}{r}
\chi_{12}\\
0\\0\\0\\0\\0\\0\\0\\0\\0\\1\\0\\1\\1\\2\\2\\5\\4\\8\\10\\16\\19\\31\\36\\57\\72\\105\\133\\195\\247\\351\\456\\633\\825\\1140\\1484\\2027\\2660\\3595\\4719\\6358\\8343\\11174\\14696\\19592\\25765\\34268\\45046\\59753\\78593\\104015\\136788\\180819\\
\end{array}$
$\begin{array}{r}
\chi_{13}\\
0\\0\\0\\0\\0\\0\\0\\0\\0\\0\\0\\1\\1\\1\\1\\3\\3\\6\\7\\11\\14\\23\\28\\44\\57\\83\\108\\158\\206\\293\\388\\541\\718\\996\\1318\\1812\\2410\\3278\\4359\\5911\\7854\\10595\\14088\\18924\\25151\\33705\\44762\\59844\\79476\\106029\\140765\\187582\\
\end{array}$
$\begin{array}{r}
\chi_{14}\\
0\\0\\0\\0\\0\\0\\0\\0\\0\\0\\0\\0\\1\\0\\1\\1\\2\\2\\5\\4\\9\\10\\17\\20\\35\\40\\65\\81\\122\\154\\231\\294\\427\\556\\786\\1031\\1448\\1900\\2635\\3487\\4784\\6341\\8659\\11484\\15592\\20729\\28007\\37257\\50210\\66802\\89792\\119568\\
\end{array}$
$\begin{array}{r}
\chi_{15}\\
0\\0\\0\\0\\0\\0\\0\\0\\0\\0\\1\\1\\1\\3\\4\\5\\8\\12\\16\\24\\34\\47\\68\\94\\129\\182\\251\\343\\474\\648\\882\\1205\\1636\\2214\\3003\\4054\\5462\\7364\\9901\\13286\\17833\\23889\\31965\\42754\\57117\\76233\\101719\\135595\\180651\\240598\\320270\\426126\\
\end{array}$
$\begin{array}{r}
\chi_{18}\\
0\\0\\0\\0\\0\\0\\0\\0\\0\\0\\0\\1\\1\\2\\2\\4\\6\\10\\13\\21\\28\\43\\59\\88\\119\\174\\237\\339\\464\\654\\891\\1247\\1697\\2350\\3195\\4398\\5964\\8166\\11055\\15066\\20367\\27660\\37321\\50545\\68110\\92012\\123852\\167013\\224568\\302400\\406322\\546554\\
\end{array}$
$\begin{array}{r}
\chi_{19}\\
0\\0\\0\\0\\0\\0\\0\\0\\0\\0\\0\\0\\1\\1\\2\\3\\5\\6\\11\\15\\24\\33\\50\\69\\104\\142\\207\\287\\411\\565\\801\\1101\\1544\\2118\\2945\\4031\\5575\\7609\\10465\\14264\\19535\\26564\\36266\\49237\\67031\\90879\\123444\\167163\\226702\\306694\\415386\\561634\\
\end{array}$
$\begin{array}{r}
\chi_{20}\\
0\\0\\0\\0\\0\\0\\0\\0\\0\\0\\0\\0\\1\\0\\1\\2\\4\\4\\8\\10\\18\\23\\37\\49\\78\\103\\156\\213\\312\\424\\614\\837\\1193\\1630\\2292\\3133\\4378\\5970\\8279\\11296\\15577\\21216\\29147\\39652\\54283\\73787\\100730\\136801\\186408\\252958\\344151\\466865\\
\end{array}$
$\begin{array}{r}
\chi_{21}\\
0\\0\\0\\0\\0\\0\\0\\0\\0\\0\\0\\0\\0\\0\\1\\0\\2\\1\\4\\4\\10\\10\\21\\25\\45\\56\\94\\122\\195\\258\\394\\532\\794\\1076\\1574\\2151\\3095\\4242\\6037\\8285\\11693\\16060\\22502\\30925\\43117\\59230\\82260\\113005\\156467\\214924\\296968\\407915\\
\end{array}$
\newpage
$\begin{array}{r}
\chi_{22}\\
0\\0\\0\\0\\0\\0\\0\\0\\0\\0\\0\\0\\0\\1\\2\\2\\4\\6\\10\\15\\23\\34\\52\\75\\110\\160\\233\\330\\475\\674\\957\\1348\\1900\\2659\\3729\\5192\\7232\\10041\\13929\\19252\\26613\\36689\\50550\\69516\\95546\\131117\\179873\\246441\\337531\\461928\\631998\\864160\\
\end{array}$
$\begin{array}{r}
\chi_{23}\\
0\\0\\0\\0\\0\\0\\0\\0\\0\\0\\0\\0\\0\\1\\1\\3\\3\\7\\9\\16\\21\\37\\49\\79\\109\\167\\231\\347\\479\\704\\978\\1409\\1953\\2789\\3858\\5451\\7536\\10565\\14577\\20324\\27982\\38826\\53386\\73781\\101287\\139598\\191387\\263152\\360459\\494762\\677230\\928510\\
\end{array}$
$\begin{array}{r}
\chi_{24}\\
0\\0\\0\\0\\0\\0\\0\\0\\0\\0\\0\\0\\0\\0\\1\\1\\2\\4\\7\\11\\18\\27\\43\\63\\98\\144\\218\\315\\467\\674\\985\\1408\\2036\\2897\\4150\\5874\\8353\\11775\\16650\\23380\\32897\\46057\\64561\\90141\\125972\\175522\\244703\\340398\\473727\\658238\\914987\\1270500\\
\end{array}$
$\begin{array}{r}
\chi_{25}\\
0\\0\\0\\0\\0\\0\\0\\0\\0\\0\\0\\0\\0\\0\\0\\1\\1\\2\\3\\5\\8\\14\\20\\33\\49\\75\\109\\168\\241\\360\\519\\761\\1090\\1584\\2254\\3244\\4607\\6573\\9292\\13198\\18589\\26267\\36912\\51956\\72836\\102248\\143060\\200386\\280048\\391639\\546869\\764184\\
\end{array}$
$\begin{array}{r}
\chi_{28}\\
0\\0\\0\\0\\0\\0\\0\\0\\0\\0\\0\\0\\0\\0\\0\\1\\2\\3\\5\\8\\14\\21\\34\\52\\80\\119\\180\\265\\394\\573\\840\\1213\\1760\\2522\\3627\\5169\\7381\\10464\\14855\\20974\\29636\\41696\\58686\\82334\\115531\\161707\\226369\\316299\\442002\\616839\\860993\\1200744\\
\end{array}$
$\begin{array}{r}
\chi_{29}\\
0\\0\\0\\0\\0\\0\\0\\0\\0\\0\\0\\0\\0\\0\\0\\0\\1\\1\\3\\4\\8\\11\\21\\30\\50\\73\\117\\169\\264\\382\\578\\836\\1244\\1789\\2631\\3769\\5479\\7827\\11286\\16062\\23021\\32670\\46576\\65946\\93657\\132322\\187394\\264367\\373589\\526554\\743134\\1046943\\
\end{array}$
$\begin{array}{r}
\chi_{30}\\
0\\0\\0\\0\\0\\0\\0\\0\\0\\0\\0\\0\\0\\0\\0\\1\\1\\3\\4\\9\\12\\23\\32\\55\\79\\128\\183\\287\\412\\627\\899\\1341\\1918\\2819\\4019\\5839\\8300\\11949\\16932\\24211\\34211\\48642\\68569\\97075\\136559\\192690\\270618\\380910\\534341\\750825\\1052524\\1477500\\
\end{array}$
$\begin{array}{r}
\chi_{31}\\
0\\0\\0\\0\\0\\0\\0\\0\\0\\0\\0\\0\\0\\0\\0\\0\\1\\2\\3\\6\\10\\16\\27\\44\\67\\107\\163\\249\\376\\567\\838\\1251\\1836\\2699\\3936\\5741\\8302\\12028\\17308\\24908\\35694\\51136\\72979\\104174\\148262\\210974\\299632\\425496\\603326\\855621\\1212253\\1718020\\
\end{array}$
$\begin{array}{r}
\chi_{32}\\
0\\0\\0\\0\\0\\0\\0\\0\\0\\0\\0\\0\\0\\0\\0\\0\\0\\0\\1\\0\\2\\2\\5\\6\\13\\17\\32\\44\\74\\105\\172\\244\\384\\554\\845\\1219\\1833\\2642\\3921\\5649\\8290\\11932\\17398\\24987\\36228\\51995\\75076\\107640\\155043\\222201\\319590\\458190\\
\end{array}$

\newpage
$\begin{array}{r}
\chi_{33}\\
0\\0\\0\\0\\0\\0\\0\\0\\0\\0\\0\\0\\0\\0\\0\\0\\1\\1\\3\\4\\9\\12\\24\\34\\59\\86\\141\\205\\326\\473\\728\\1058\\1595\\2308\\3431\\4949\\7265\\10448\\15202\\21785\\31486\\44994\\64663\\92196\\131956\\187758\\267927\\380704\\542085\\769649\\1094568\\1553736\\
\end{array}$
$\begin{array}{r}
\chi_{34}\\
0\\0\\0\\0\\0\\0\\0\\0\\0\\0\\0\\0\\0\\0\\0\\0\\0\\0\\0\\0\\1\\0\\2\\2\\6\\6\\15\\18\\38\\50\\91\\126\\217\\304\\499\\717\\1134\\1639\\2537\\3683\\5609\\8162\\12269\\17882\\26660\\38858\\57578\\84001\\123965\\180992\\266561\\389774\\
\end{array}$
$\begin{array}{r}
\chi_{35}\\
0\\0\\0\\0\\0\\0\\0\\0\\0\\0\\0\\0\\0\\0\\0\\0\\0\\0\\0\\0\\0\\0\\0\\0\\0\\0\\0\\0\\0\\0\\1\\0\\3\\2\\7\\9\\23\\30\\65\\95\\179\\272\\487\\751\\1291\\2020\\3367\\5308\\8719\\13777\\22420\\35613\\
\end{array}$
$\begin{array}{r}
\chi_{36}\\
0\\0\\0\\0\\0\\0\\0\\0\\0\\0\\0\\0\\0\\0\\0\\0\\0\\0\\2\\1\\4\\6\\12\\16\\33\\45\\80\\117\\191\\279\\449\\650\\1009\\1478\\2239\\3262\\4890\\7096\\10505\\15232\\22335\\32286\\47083\\67875\\98473\\141781\\204934\\294657\\425017\\610691\\879754\\1264463\\
\end{array}$
$\begin{array}{r}
\chi_{37}\\
0\\0\\0\\0\\0\\0\\0\\0\\0\\0\\0\\0\\0\\0\\0\\0\\0\\1\\1\\4\\6\\12\\18\\34\\51\\87\\131\\215\\323\\509\\758\\1171\\1735\\2625\\3866\\5779\\8464\\12514\\18242\\26759\\38850\\56628\\81924\\118856\\171494\\247896\\356965\\514746\\740288\\1065918\\1532269\\2205293\\
\end{array}$
$\begin{array}{r}
\chi_{38}\\
0\\0\\0\\0\\0\\0\\0\\0\\0\\0\\0\\0\\0\\0\\0\\0\\0\\1\\2\\4\\7\\14\\23\\39\\62\\104\\163\\258\\396\\616\\935\\1425\\2136\\3214\\4774\\7099\\10464\\15433\\22607\\33101\\48233\\70255\\101949\\147854\\213856\\309191\\446165\\643632\\927366\\1336177\\1924034\\2771428\\
\end{array}$
$\begin{array}{r}
\chi_{39}\\
0\\0\\0\\0\\0\\0\\0\\0\\0\\0\\0\\0\\0\\0\\0\\0\\0\\0\\0\\0\\0\\0\\0\\0\\0\\0\\0\\0\\0\\0\\0\\0\\0\\0\\0\\0\\0\\0\\0\\0\\0\\0\\0\\0\\0\\2\\11\\33\\98\\269\\685\\1676\\
\end{array}$
$\begin{array}{r}
\chi_{41}\\
0\\0\\0\\0\\0\\0\\0\\0\\0\\0\\0\\0\\0\\0\\0\\0\\0\\0\\0\\0\\0\\0\\0\\0\\0\\0\\0\\0\\0\\0\\0\\0\\0\\0\\0\\0\\0\\0\\0\\0\\0\\0\\0\\0\\0\\0\\2\\10\\47\\147\\477\\1302\\
\end{array}$
$\begin{array}{r}
\chi_{43}\\
0\\0\\0\\0\\0\\0\\0\\0\\0\\0\\0\\0\\0\\0\\0\\0\\0\\0\\1\\2\\3\\7\\11\\21\\35\\59\\94\\155\\242\\384\\594\\921\\1404\\2149\\3237\\4888\\7311\\10929\\16224\\24084\\35541\\52434\\77054\\113157\\165707\\242576\\354369\\517625\\755213\\1102014\\1607414\\2346206\\
\end{array}$
\newpage
$\begin{array}{r}
\chi_{44}\\
0\\0\\0\\0\\0\\0\\0\\0\\0\\0\\0\\0\\0\\0\\0\\0\\0\\0\\0\\0\\0\\0\\0\\0\\0\\0\\0\\0\\0\\0\\0\\0\\0\\0\\0\\0\\0\\0\\0\\0\\0\\0\\0\\0\\0\\0\\6\\24\\107\\330\\1025\\2812\\
\end{array}$
$\begin{array}{r}
\chi_{46}\\
0\\0\\0\\0\\0\\0\\0\\0\\0\\0\\0\\0\\0\\0\\0\\0\\0\\0\\1\\1\\3\\7\\12\\21\\38\\60\\102\\165\\265\\416\\660\\1015\\1574\\2405\\3665\\5535\\8364\\12512\\18726\\27862\\41391\\61240\\90541\\133381\\196395\\288534\\423656\\621219\\910885\\1334634\\1956465\\2868528\\
\end{array}$
$\begin{array}{r}
\chi_{47}\\
0\\0\\0\\0\\0\\0\\0\\0\\0\\0\\0\\0\\0\\0\\0\\0\\0\\0\\0\\0\\0\\0\\0\\0\\0\\0\\0\\0\\0\\0\\0\\0\\0\\0\\0\\0\\0\\0\\0\\0\\1\\5\\9\\32\\71\\166\\360\\810\\1679\\3593\\7392\\15386\\
\end{array}$
$\begin{array}{r}
\chi_{49}\\
0\\0\\0\\0\\0\\0\\0\\0\\0\\0\\0\\0\\0\\0\\0\\0\\0\\0\\0\\0\\1\\1\\2\\3\\8\\10\\21\\33\\57\\87\\149\\224\\366\\559\\881\\1337\\2083\\3134\\4810\\7231\\10967\\16415\\24746\\36907\\55336\\82412\\123111\\183173\\273314\\406796\\607295\\906291\\
\end{array}$
$\begin{array}{r}
\chi_{50}\\
0\\0\\0\\0\\0\\0\\0\\0\\0\\0\\0\\0\\0\\0\\0\\0\\0\\0\\0\\0\\0\\1\\0\\2\\2\\5\\7\\17\\21\\43\\63\\111\\165\\283\\419\\688\\1040\\1651\\2494\\3909\\5887\\9104\\13735\\21026\\31694\\48303\\72780\\110592\\166965\\253636\\384221\\585754\\
\end{array}$
$\begin{array}{r}
\chi_{51}\\
0\\0\\0\\0\\0\\0\\0\\0\\0\\0\\0\\0\\0\\0\\0\\0\\0\\0\\0\\0\\0\\0\\0\\0\\0\\0\\0\\0\\0\\0\\0\\0\\0\\0\\0\\0\\2\\2\\6\\16\\36\\68\\147\\283\\563\\1074\\2054\\3863\\7301\\13626\\25591\\48083\\
\end{array}$
$\begin{array}{r}
\chi_{52}\\
0\\0\\0\\0\\0\\0\\0\\0\\0\\0\\0\\0\\0\\0\\0\\0\\0\\0\\0\\2\\4\\7\\13\\25\\43\\75\\123\\205\\329\\530\\831\\1312\\2030\\3139\\4800\\7334\\11096\\16780\\25191\\37778\\56370\\84006\\124713\\185003\\273668\\404562\\597013\\880792\\1298177\\1913903\\2821306\\4162826\\
\end{array}$
$\begin{array}{r}
\chi_{53}\\
0\\0\\0\\0\\0\\0\\0\\0\\0\\0\\0\\0\\0\\0\\0\\0\\0\\0\\0\\0\\0\\0\\0\\0\\0\\0\\0\\0\\0\\0\\0\\0\\0\\0\\0\\0\\0\\0\\0\\0\\0\\0\\0\\0\\1\\9\\41\\140\\429\\1248\\3371\\8771\\
\end{array}$
$\begin{array}{r}
\chi_{55}\\
0\\0\\0\\0\\0\\0\\0\\0\\0\\0\\0\\0\\0\\0\\0\\0\\0\\0\\0\\0\\0\\0\\0\\0\\0\\0\\0\\0\\0\\0\\0\\0\\0\\0\\0\\0\\0\\0\\0\\0\\0\\0\\0\\0\\0\\0\\10\\47\\198\\678\\2160\\6233\\
\end{array}$
$\begin{array}{r}
\chi_{57}\\
0\\0\\0\\0\\0\\0\\0\\0\\0\\0\\0\\0\\0\\0\\0\\0\\0\\0\\0\\0\\1\\2\\5\\8\\17\\28\\52\\85\\147\\236\\393\\620\\1000\\1562\\2468\\3811\\5928\\9077\\13954\\21215\\32327\\48873\\73990\\111386\\167844\\251982\\378615\\567613\\851926\\1277559\\1919182\\2884882\\
\end{array}$
\newpage
$\begin{array}{r}
\chi_{58}\\
0\\0\\0\\0\\0\\0\\0\\0\\0\\0\\0\\0\\0\\0\\0\\0\\0\\0\\0\\1\\2\\4\\9\\16\\29\\51\\87\\143\\239\\384\\617\\977\\1539\\2391\\3715\\5705\\8742\\13306\\20201\\30494\\45970\\68989\\103391\\154497\\230610\\343506\\511443\\760563\\1131101\\1681860\\2502839\\3727955\\
\end{array}$
$\begin{array}{r}
\chi_{59}\\
0\\0\\0\\0\\0\\0\\0\\0\\0\\0\\0\\0\\0\\0\\0\\0\\0\\0\\0\\0\\0\\0\\0\\0\\0\\0\\0\\0\\0\\0\\0\\0\\0\\0\\0\\0\\0\\0\\0\\0\\0\\0\\0\\1\\5\\23\\85\\291\\903\\2689\\7533\\20410\\
\end{array}$
$\begin{array}{r}
\chi_{61}\\
0\\0\\0\\0\\0\\0\\0\\0\\0\\0\\0\\0\\0\\0\\0\\0\\0\\0\\0\\0\\0\\1\\2\\4\\9\\16\\30\\54\\93\\159\\267\\440\\718\\1162\\1858\\2952\\4661\\7306\\11400\\17710\\27390\\42228\\64923\\99552\\152385\\232927\\355691\\542958\\829017\\1266777\\1938934\\2975374\\
\end{array}$
$\begin{array}{r}
\chi_{62}\\
0\\0\\0\\0\\0\\0\\0\\0\\0\\0\\0\\0\\0\\0\\0\\0\\0\\0\\0\\0\\1\\2\\6\\10\\21\\36\\68\\112\\197\\320\\537\\857\\1395\\2196\\3499\\5445\\8530\\13158\\20367\\31172\\47816\\72766\\110862\\167984\\254744\\384927\\582147\\878626\\1327792\\2005663\\3036098\\4602110\\
\end{array}$
$\begin{array}{r}
\chi_{63}\\
0\\0\\0\\0\\0\\0\\0\\0\\0\\0\\0\\0\\0\\0\\0\\0\\0\\0\\0\\0\\1\\1\\4\\8\\16\\29\\54\\90\\160\\263\\440\\713\\1164\\1845\\2952\\4628\\7272\\11289\\17537\\26987\\41557\\63560\\97197\\148001\\225331\\342055\\519468\\787727\\1195664\\1815025\\2760934\\4207860\\
\end{array}$
$\begin{array}{r}
\chi_{64}\\
0\\0\\0\\0\\0\\0\\0\\0\\0\\0\\0\\0\\0\\0\\0\\0\\0\\0\\0\\0\\0\\1\\2\\4\\9\\16\\30\\54\\92\\158\\266\\437\\714\\1157\\1851\\2943\\4652\\7297\\11397\\17724\\27437\\42351\\65190\\100086\\153415\\234871\\359272\\549546\\841133\\1289231\\1981167\\3056268\\
\end{array}$
$\begin{array}{r}
\chi_{65}\\
0\\0\\0\\0\\0\\0\\0\\0\\0\\0\\0\\0\\0\\0\\0\\0\\0\\0\\0\\0\\0\\1\\1\\5\\8\\19\\30\\63\\100\\185\\299\\519\\829\\1391\\2198\\3585\\5633\\9000\\14033\\22134\\34300\\53546\\82638\\128114\\197128\\304405\\467848\\721435\\1110256\\1715028\\2651226\\4120694\\
\end{array}$
$\begin{array}{r}
\chi_{66}\\
0\\0\\0\\0\\0\\0\\0\\0\\0\\0\\0\\0\\0\\0\\0\\0\\0\\0\\0\\0\\0\\0\\1\\1\\5\\6\\16\\25\\53\\85\\160\\257\\455\\730\\1238\\1976\\3262\\5168\\8366\\13191\\21062\\33047\\52299\\81810\\128755\\201224\\316045\\494757\\778647\\1225773\\1943312\\3095438\\
\end{array}$
\newpage
$\begin{array}{r}
\chi_{67}\\
0\\0\\0\\0\\0\\0\\0\\0\\0\\0\\0\\0\\0\\0\\0\\0\\0\\0\\0\\0\\0\\2\\4\\10\\17\\37\\64\\119\\201\\355\\582\\984\\1588\\2603\\4146\\6656\\10469\\16558\\25812\\40337\\62437\\96771\\148967\\229545\\352123\\540664\\828112\\1270153\\1946953\\2991932\\4605230\\7120281\\
\end{array}$
$\begin{array}{r}
\chi_{68}\\
0\\0\\0\\0\\0\\0\\0\\0\\0\\0\\0\\0\\0\\0\\0\\0\\0\\0\\0\\0\\0\\0\\1\\2\\5\\10\\22\\38\\74\\127\\228\\381\\657\\1077\\1800\\2910\\4753\\7599\\12222\\19354\\30771\\48415\\76353\\119590\\187678\\293308\\459417\\718418\\1127044\\1770153\\2794895\\4434008\\
\end{array}$
$\begin{array}{r}
\chi_{69}\\
0\\0\\0\\0\\0\\0\\0\\0\\0\\0\\0\\0\\0\\0\\0\\0\\0\\0\\0\\0\\0\\0\\0\\0\\0\\0\\0\\0\\0\\0\\0\\0\\0\\0\\0\\0\\1\\0\\3\\6\\18\\35\\101\\208\\535\\1246\\3066\\7293\\17771\\42347\\101467\\240033\\
\end{array}$
$\begin{array}{r}
\chi_{70}\\
0\\0\\0\\0\\0\\0\\0\\0\\0\\0\\0\\0\\0\\0\\0\\0\\0\\0\\0\\0\\0\\0\\0\\0\\0\\0\\0\\0\\0\\0\\0\\0\\0\\0\\0\\0\\4\\6\\24\\55\\144\\306\\722\\1508\\3299\\6844\\14472\\29909\\62712\\130290\\273986\\576333\\
\end{array}$
$\begin{array}{r}
\chi_{71}\\
0\\0\\0\\0\\0\\0\\0\\0\\0\\0\\0\\0\\0\\0\\0\\0\\0\\0\\0\\0\\0\\0\\0\\0\\0\\0\\0\\0\\0\\0\\0\\0\\0\\0\\0\\0\\0\\0\\0\\4\\14\\50\\140\\400\\1018\\2578\\6203\\14857\\34638\\80399\\184107\\420235\\
\end{array}$
$\begin{array}{r}
\chi_{73}\\
0\\0\\0\\0\\0\\0\\0\\0\\0\\0\\0\\0\\0\\0\\0\\0\\0\\0\\0\\0\\0\\0\\0\\0\\0\\0\\0\\0\\0\\0\\3\\2\\10\\19\\46\\81\\179\\317\\636\\1139\\2153\\3841\\7090\\12575\\22835\\40668\\73558\\132129\\241081\\440946\\820710\\1543715\\
\end{array}$
$\begin{array}{r}
\chi_{74}\\
0\\0\\0\\0\\0\\0\\0\\0\\0\\0\\0\\0\\0\\0\\0\\0\\0\\0\\0\\0\\0\\0\\0\\0\\0\\0\\0\\0\\0\\0\\0\\0\\0\\0\\0\\0\\0\\0\\0\\0\\1\\13\\46\\182\\552\\1672\\4602\\12476\\32209\\82015\\203055\\496080\\
\end{array}$
$\begin{array}{r}
\chi_{76}\\
0\\0\\0\\0\\0\\0\\0\\0\\0\\0\\0\\0\\0\\0\\0\\0\\0\\0\\0\\0\\0\\0\\1\\3\\6\\15\\29\\57\\106\\194\\339\\595\\1012\\1712\\2849\\4715\\7700\\12539\\20223\\32510\\51950\\82808\\131436\\208390\\329621\\521393\\824624\\1306649\\2075513\\3312772\\5321185\\8627638\\
\end{array}$
$\begin{array}{r}
\chi_{77}\\
0\\0\\0\\0\\0\\0\\0\\0\\0\\0\\0\\0\\0\\0\\0\\0\\0\\0\\0\\0\\0\\0\\1\\2\\6\\13\\28\\51\\99\\176\\319\\548\\951\\1594\\2685\\4419\\7282\\11826\\19194\\30831\\49498\\78920\\125770\\199592\\316823\\501973\\796597\\1265056\\2016658\\3228408\\5207004\\8476067\\
\end{array}$
\newpage
$\begin{array}{r}
\chi_{78}\\
0\\0\\0\\0\\0\\0\\0\\0\\0\\0\\0\\0\\0\\0\\0\\0\\0\\0\\0\\0\\0\\0\\0\\1\\2\\6\\11\\26\\46\\92\\161\\296\\505\\889\\1488\\2534\\4184\\6964\\11363\\18632\\30138\\48895\\78653\\126807\\203504\\327559\\526721\\851211\\1380465\\2258698\\3731422\\6257769\\
\end{array}$
$\begin{array}{r}
\chi_{79}\\
0\\0\\0\\0\\0\\0\\0\\0\\0\\0\\0\\0\\0\\0\\0\\0\\0\\0\\0\\0\\0\\0\\0\\0\\2\\3\\9\\17\\38\\66\\133\\230\\425\\724\\1275\\2138\\3659\\6052\\10132\\16613\\27413\\44625\\72993\\118392\\192814\\312805\\509922\\831670\\1366442\\2257782\\3774641\\6388609\\
\end{array}$
$\begin{array}{r}
\chi_{80}\\
0\\0\\0\\0\\0\\0\\0\\0\\0\\0\\0\\0\\0\\0\\0\\0\\0\\0\\0\\0\\0\\0\\0\\2\\5\\11\\23\\48\\91\\170\\305\\543\\941\\1617\\2727\\4574\\7576\\12469\\20361\\33104\\53498\\86200\\138383\\221791\\354895\\567930\\909450\\1460102\\2353200\\3816367\\6242657\\10333106\\
\end{array}$
$\begin{array}{r}
\chi_{81}\\
0\\0\\0\\0\\0\\0\\0\\0\\0\\0\\0\\0\\0\\0\\0\\0\\0\\0\\0\\0\\0\\0\\0\\0\\0\\0\\0\\0\\0\\0\\0\\0\\0\\0\\0\\0\\0\\0\\0\\1\\8\\46\\150\\521\\1500\\4283\\11410\\29945\\75632\\188680\\460125\\1109339\\
\end{array}$
$\begin{array}{r}
\chi_{83}\\
0\\0\\0\\0\\0\\0\\0\\0\\0\\0\\0\\0\\0\\0\\0\\0\\0\\0\\0\\0\\0\\0\\0\\0\\0\\0\\0\\0\\0\\0\\0\\0\\0\\0\\0\\0\\0\\0\\0\\0\\1\\7\\46\\191\\717\\2333\\7190\\20665\\57174\\151862\\392584\\988361\\
\end{array}$
$\begin{array}{r}
\chi_{85}\\
0\\0\\0\\0\\0\\0\\0\\0\\0\\0\\0\\0\\0\\0\\0\\0\\0\\0\\0\\0\\0\\0\\0\\0\\0\\0\\0\\0\\0\\0\\0\\0\\0\\0\\0\\0\\0\\0\\0\\0\\2\\18\\75\\313\\1039\\3298\\9652\\27194\\73200\\191843\\488400\\1218333\\
\end{array}$
$\begin{array}{r}
\chi_{87}\\
0\\0\\0\\0\\0\\0\\0\\0\\0\\0\\0\\0\\0\\0\\0\\0\\0\\0\\0\\0\\0\\0\\0\\0\\0\\0\\0\\0\\0\\0\\0\\0\\0\\0\\0\\0\\3\\4\\21\\50\\151\\358\\951\\2222\\5470\\12725\\30150\\69636\\162141\\372591\\858667\\1964610\\
\end{array}$
$\begin{array}{r}
\chi_{88}\\
0\\0\\0\\0\\0\\0\\0\\0\\0\\0\\0\\0\\0\\0\\0\\0\\0\\0\\0\\0\\0\\0\\0\\0\\0\\0\\0\\0\\0\\0\\0\\0\\0\\0\\0\\0\\0\\0\\2\\13\\51\\171\\506\\1416\\3744\\9597\\23903\\58545\\141005\\335833\\791626\\1851089\\
\end{array}$
\newpage
$\begin{array}{r}
\chi_{89}\\
0\\0\\0\\0\\0\\0\\0\\0\\0\\0\\0\\0\\0\\0\\0\\0\\0\\0\\0\\0\\0\\0\\0\\0\\0\\0\\0\\0\\0\\0\\0\\0\\0\\0\\0\\0\\0\\3\\11\\43\\121\\353\\919\\2353\\5784\\14077\\33499\\79245\\185407\\431815\\999294\\2303725\\
\end{array}$
$\begin{array}{r}
\chi_{91}\\
0\\0\\0\\0\\0\\0\\0\\0\\0\\0\\0\\0\\0\\0\\0\\0\\0\\0\\0\\0\\0\\0\\0\\0\\0\\0\\0\\0\\0\\0\\1\\1\\4\\9\\19\\40\\85\\161\\320\\614\\1170\\2223\\4250\\8124\\15771\\31057\\62227\\127273\\265407\\563819\\1216596\\2659383\\
\end{array}$
$\begin{array}{r}
\chi_{92}\\
0\\0\\0\\0\\0\\0\\0\\0\\0\\0\\0\\0\\0\\0\\0\\0\\0\\0\\0\\0\\0\\0\\0\\0\\1\\3\\7\\15\\31\\61\\115\\214\\384\\683\\1190\\2055\\3495\\5919\\9909\\16531\\27380\\45257\\74484\\122607\\201608\\332510\\550246\\917414\\1544048\\2635887\\4578224\\8128645\\
\end{array}$
$\begin{array}{r}
\chi_{93}\\
0\\0\\0\\0\\0\\0\\0\\0\\0\\0\\0\\0\\0\\0\\0\\0\\0\\0\\0\\0\\0\\0\\0\\1\\3\\8\\17\\38\\72\\141\\259\\473\\832\\1466\\2507\\4285\\7199\\12054\\19952\\32957\\53990\\88340\\143861\\234222\\380675\\619906\\1011030\\1657465\\2734112\\4555888\\7687351\\13192303\\
\end{array}$
$\begin{array}{r}
\chi_{94}\\
0\\0\\0\\0\\0\\0\\0\\0\\0\\0\\0\\0\\0\\0\\0\\0\\0\\0\\0\\0\\0\\0\\0\\0\\0\\1\\3\\7\\17\\35\\70\\140\\264\\486\\887\\1583\\2794\\4880\\8447\\14494\\24777\\42098\\71381\\120793\\204521\\346712\\590446\\1011336\\1748606\\3059265\\5435769\\9835083\\
\end{array}$
$\begin{array}{r}
\chi_{95}\\
0\\0\\0\\0\\0\\0\\0\\0\\0\\0\\0\\0\\0\\0\\0\\0\\0\\0\\0\\0\\0\\0\\0\\0\\1\\2\\5\\12\\24\\48\\94\\174\\320\\571\\1009\\1753\\3028\\5149\\8728\\14664\\24561\\40912\\68123\\113168\\188444\\314483\\528165\\894243\\1533379\\2670540\\4746410\\8633692\\
\end{array}$
$\begin{array}{r}
\chi_{96}\\
0\\0\\0\\0\\0\\0\\0\\0\\0\\0\\0\\0\\0\\0\\0\\0\\0\\0\\0\\0\\0\\0\\0\\0\\0\\0\\0\\0\\1\\1\\3\\5\\13\\20\\44\\76\\147\\253\\472\\810\\1482\\2579\\4715\\8510\\16147\\30994\\62711\\130101\\281192\\620304\\1399545\\3185945\\
\end{array}$
$\begin{array}{r}
\chi_{97}\\
0\\0\\0\\0\\0\\0\\0\\0\\0\\0\\0\\0\\0\\0\\0\\0\\0\\0\\0\\0\\0\\0\\0\\0\\2\\4\\11\\24\\50\\96\\189\\342\\626\\1106\\1943\\3340\\5731\\9668\\16272\\27121\\45075\\74477\\122944\\202346\\333401\\549786\\910390\\1515705\\2546756\\4329221\\7477233\\13163453\\
\end{array}$
\newpage
$\begin{array}{r}
\chi_{98}\\
0\\0\\0\\0\\0\\0\\0\\0\\0\\0\\0\\0\\0\\0\\0\\0\\0\\0\\0\\0\\0\\0\\0\\0\\0\\0\\0\\0\\0\\0\\4\\3\\16\\28\\74\\133\\298\\543\\1113\\2039\\3967\\7281\\13859\\25554\\48409\\90717\\173977\\334919\\659976\\1316915\\2686648\\5563486\\
\end{array}$
$\begin{array}{r}
\chi_{99}\\
0\\0\\0\\0\\0\\0\\0\\0\\0\\0\\0\\0\\0\\0\\0\\0\\0\\0\\0\\0\\0\\0\\0\\0\\0\\0\\0\\0\\0\\0\\0\\0\\0\\0\\0\\0\\0\\0\\1\\9\\46\\175\\604\\1886\\5591\\15633\\42402\\111122\\284770\\713525\\1758319\\4263252\\
\end{array}$
$\begin{array}{r}
\chi_{101}\\
0\\0\\0\\0\\0\\0\\0\\0\\0\\0\\0\\0\\0\\0\\0\\0\\0\\0\\0\\0\\0\\0\\0\\0\\1\\1\\6\\11\\29\\51\\110\\196\\382\\667\\1225\\2111\\3734\\6348\\10953\\18446\\31343\\52484\\88478\\148095\\249678\\420811\\716699\\1229757\\2144617\\3801040\\6893845\\12801184\\
\end{array}$
$\begin{array}{r}
\chi_{102}\\
0\\0\\0\\0\\0\\0\\0\\0\\0\\0\\0\\0\\0\\0\\0\\0\\0\\0\\0\\0\\0\\0\\0\\0\\0\\0\\0\\0\\0\\0\\0\\0\\0\\0\\0\\0\\0\\0\\3\\21\\85\\296\\936\\2723\\7563\\20260\\52637\\133710\\333338\\817725\\1978892\\4731646\\
\end{array}$
$\begin{array}{r}
\chi_{104}\\
0\\0\\0\\0\\0\\0\\0\\0\\0\\0\\0\\0\\0\\0\\0\\0\\0\\0\\0\\0\\0\\0\\0\\0\\0\\1\\4\\6\\17\\31\\66\\122\\238\\421\\777\\1355\\2400\\4127\\7151\\12138\\20740\\35021\\59434\\100499\\171215\\292759\\507251\\890248\\1594972\\2921759\\5499070\\10634085\\
\end{array}$
$\begin{array}{r}
\chi_{105}\\
0\\0\\0\\0\\0\\0\\0\\0\\0\\0\\0\\0\\0\\0\\0\\0\\0\\0\\0\\0\\0\\0\\0\\0\\0\\0\\0\\0\\0\\0\\0\\0\\0\\0\\0\\0\\0\\0\\0\\3\\26\\136\\514\\1761\\5427\\15920\\44178\\118700\\308656\\784698\\1951765\\4773670\\
\end{array}$
$\begin{array}{r}
\chi_{107}\\
0\\0\\0\\0\\0\\0\\0\\0\\0\\0\\0\\0\\0\\0\\0\\0\\0\\0\\0\\0\\0\\0\\0\\0\\0\\0\\0\\0\\0\\0\\0\\0\\0\\0\\0\\0\\0\\0\\0\\3\\26\\136\\514\\1761\\5427\\15920\\44178\\118700\\308656\\784698\\1951765\\4773669\\
\end{array}$
$\begin{array}{r}
\chi_{109}\\
0\\0\\0\\0\\0\\0\\0\\0\\0\\0\\0\\0\\0\\0\\0\\0\\0\\0\\0\\0\\0\\0\\0\\0\\0\\0\\0\\0\\0\\0\\5\\6\\25\\48\\121\\233\\515\\983\\2010\\3806\\7460\\14066\\27098\\51225\\98594\\189119\\369296\\726398\\1456587\\2959329\\6124162\\12848040\\
\end{array}$
\newpage
$\begin{array}{r}
\chi_{110}\\
0\\0\\0\\0\\0\\0\\0\\0\\0\\0\\0\\0\\0\\0\\0\\0\\0\\0\\0\\0\\0\\0\\0\\0\\0\\1\\5\\9\\26\\51\\112\\214\\425\\778\\1460\\2606\\4701\\8241\\14516\\25099\\43569\\74836\\129038\\221635\\383114\\664032\\1163714\\2060522\\3712178\\6810236\\12785353\\24565376\\
\end{array}$
$\begin{array}{r}
\chi_{111}\\
0\\0\\0\\0\\0\\0\\0\\0\\0\\0\\0\\0\\0\\0\\0\\0\\0\\0\\0\\0\\0\\0\\0\\0\\0\\2\\5\\15\\31\\74\\144\\295\\551\\1054\\1904\\3481\\6144\\10905\\18939\\32977\\56661\\97610\\166986\\286817\\492247\\851162\\1480599\\2609337\\4662599\\8497213\\15814298\\30158915\\
\end{array}$
$\begin{array}{r}
\chi_{112}\\
0\\0\\0\\0\\0\\0\\0\\0\\0\\0\\0\\0\\0\\0\\0\\0\\0\\0\\0\\0\\0\\0\\0\\0\\0\\2\\4\\13\\28\\64\\126\\258\\484\\921\\1676\\3054\\5411\\9599\\16719\\29108\\50176\\86502\\148461\\255591\\440530\\764954\\1339198\\2377102\\4287592\\7896691\\14877219\\28736994\\
\end{array}$
$\begin{array}{r}
\chi_{113}\\
0\\0\\0\\0\\0\\0\\0\\0\\0\\0\\0\\0\\0\\0\\0\\0\\0\\0\\0\\0\\0\\0\\0\\0\\0\\0\\0\\0\\0\\0\\0\\0\\0\\0\\0\\0\\2\\6\\30\\104\\356\\1091\\3256\\9155\\25174\\66899\\174239\\443887\\1112493\\2743423\\6675835\\16039355\\
\end{array}$
$\begin{array}{r}
\chi_{114}\\
0\\0\\0\\0\\0\\0\\0\\0\\0\\0\\0\\0\\0\\0\\0\\0\\0\\0\\0\\0\\0\\0\\0\\0\\0\\0\\0\\0\\0\\0\\0\\0\\3\\6\\27\\60\\164\\362\\870\\1860\\4169\\8834\\19156\\40589\\87466\\187210\\406450\\884221\\1945417\\4299814\\9576576\\21400281\\
\end{array}$
$\begin{array}{r}
\chi_{115}\\
0\\0\\0\\0\\0\\0\\0\\0\\0\\0\\0\\0\\0\\0\\0\\0\\0\\0\\0\\0\\0\\0\\0\\0\\0\\0\\0\\0\\0\\0\\0\\1\\3\\12\\28\\82\\183\\441\\959\\2143\\4563\\9871\\20906\\44794\\95535\\206033\\445855\\975163\\2145685\\4760353\\10613730\\23777156\\
\end{array}$
$\begin{array}{r}
\chi_{116}\\
0\\0\\0\\0\\0\\0\\0\\0\\0\\0\\0\\0\\0\\0\\0\\0\\0\\0\\0\\0\\0\\0\\0\\0\\0\\0\\1\\3\\8\\20\\45\\95\\194\\383\\736\\1390\\2584\\4735\\8601\\15516\\27854\\49981\\89853\\162587\\297209\\551661\\1043446\\2019450\\4004938\\8143878\\16946351\\35982090\\
\end{array}$
$\begin{array}{r}
\chi_{117}\\
0\\0\\0\\0\\0\\0\\0\\0\\0\\0\\0\\0\\0\\0\\0\\0\\0\\0\\0\\0\\0\\0\\0\\0\\0\\0\\1\\3\\9\\22\\50\\105\\216\\425\\826\\1554\\2902\\5322\\9703\\17502\\31529\\56604\\102026\\184695\\338083\\627083\\1185268\\2287578\\4522134\\9152846\\18953912\\40031881\\
\end{array}$
\newpage
$\begin{array}{r}
\chi_{118}\\
0\\0\\0\\0\\0\\0\\0\\0\\0\\0\\0\\0\\0\\0\\0\\0\\0\\0\\0\\0\\0\\0\\0\\0\\0\\0\\1\\3\\10\\25\\57\\124\\259\\515\\1012\\1931\\3637\\6740\\12395\\22561\\40944\\73998\\133925\\242966\\444013\\819010\\1532122\\2913431\\5650041\\11186020\\22622211\\46665631\\
\end{array}$
$\begin{array}{r}
\chi_{119}\\
0\\0\\0\\0\\0\\0\\0\\0\\0\\0\\0\\0\\0\\0\\0\\0\\0\\0\\0\\0\\0\\0\\0\\0\\0\\0\\0\\0\\0\\0\\0\\0\\0\\0\\1\\2\\16\\37\\134\\373\\1109\\3005\\8269\\21672\\56604\\144111\\363358\\901132\\2212555\\5364868\\12885469\\30632762\\
\end{array}$
$\begin{array}{r}
\chi_{120}\\
0\\0\\0\\0\\0\\0\\0\\0\\0\\0\\0\\0\\0\\0\\0\\0\\0\\0\\0\\0\\0\\0\\0\\0\\0\\0\\1\\4\\12\\28\\66\\139\\291\\572\\1120\\2116\\3978\\7312\\13400\\24254\\43879\\78971\\142681\\258390\\472453\\872856\\1639253\\3133988\\6122592\\12223180\\24947456\\51933475\\
\end{array}$
$\begin{array}{r}
\chi_{121}\\
0\\0\\0\\0\\0\\0\\0\\0\\0\\0\\0\\0\\0\\0\\0\\0\\0\\0\\0\\0\\0\\0\\0\\0\\0\\0\\1\\5\\13\\34\\75\\163\\332\\668\\1286\\2461\\4581\\8481\\15453\\28080\\50585\\91232\\164306\\297772\\543111\\1003584\\1881618\\3598204\\7025278\\14035510\\28658535\\59730244\\
\end{array}$
$\begin{array}{r}
\chi_{122}\\
0\\0\\0\\0\\0\\0\\0\\0\\0\\0\\0\\0\\0\\0\\0\\0\\0\\0\\0\\0\\0\\0\\0\\0\\0\\0\\0\\0\\0\\0\\0\\0\\0\\0\\0\\1\\12\\34\\135\\414\\1278\\3608\\10145\\27114\\71712\\184544\\468863\\1170223\\2887279\\7028980\\16935961\\40366524\\
\end{array}$
$\begin{array}{r}
\chi_{123}\\
0\\0\\0\\0\\0\\0\\0\\0\\0\\0\\0\\0\\0\\0\\0\\0\\0\\0\\0\\0\\0\\0\\0\\0\\0\\0\\0\\0\\0\\0\\3\\5\\19\\42\\104\\221\\496\\1011\\2134\\4334\\8951\\18305\\38126\\79788\\170261\\367939\\808612\\1798946\\4048562\\9179756\\20926386\\47817157\\
\end{array}$
$\begin{array}{r}
\chi_{124}\\
0\\0\\0\\0\\0\\0\\0\\0\\0\\0\\0\\0\\0\\0\\0\\0\\0\\0\\0\\0\\0\\0\\0\\0\\0\\0\\0\\0\\0\\0\\0\\0\\0\\0\\0\\0\\0\\4\\27\\129\\521\\1848\\5995\\18216\\52617\\145975\\391924\\1024135\\2615616\\6550659\\16129203\\39123510\\
\end{array}$
$\begin{array}{r}
\chi_{126}\\
0\\0\\0\\0\\0\\0\\0\\0\\0\\0\\0\\0\\0\\0\\0\\0\\0\\0\\0\\0\\0\\0\\0\\0\\0\\0\\0\\0\\0\\0\\0\\0\\0\\0\\0\\4\\16\\58\\189\\589\\1700\\4802\\13119\\35123\\92160\\237937\\604598\\1515562\\3749363\\9165270\\22151465\\52974037\\
\end{array}$
\newpage
$\begin{array}{r}
\chi_{127}\\
0\\0\\0\\0\\0\\0\\0\\0\\0\\0\\0\\0\\0\\0\\0\\0\\0\\0\\0\\0\\0\\0\\0\\0\\0\\0\\0\\0\\0\\0\\4\\8\\28\\66\\161\\351\\793\\1644\\3492\\7175\\14852\\30401\\62979\\130561\\274645\\583082\\1256303\\2737697\\6039512\\13441189\\30141654\\67901554\\
\end{array}$
$\begin{array}{r}
\chi_{128}\\
0\\0\\0\\0\\0\\0\\0\\0\\0\\0\\0\\0\\0\\0\\0\\0\\0\\0\\0\\0\\0\\0\\0\\0\\0\\0\\0\\0\\0\\0\\0\\0\\0\\0\\0\\0\\0\\8\\51\\229\\878\\2951\\9280\\27396\\77623\\211748\\561897\\1453707\\3686131\\9176737\\22495571\\54369815\\
\end{array}$
$\begin{array}{r}
\chi_{130}\\
0\\0\\0\\0\\0\\0\\0\\0\\0\\0\\0\\0\\0\\0\\0\\0\\0\\0\\0\\0\\0\\0\\0\\0\\0\\0\\0\\0\\1\\2\\6\\15\\38\\81\\182\\381\\807\\1638\\3379\\6832\\14009\\28579\\59282\\123669\\262640\\563990\\1230956\\2716606\\6066205\\13646829\\30892256\\70137454\\
\end{array}$
$\begin{array}{r}
\chi_{131}\\
0\\0\\0\\0\\0\\0\\0\\0\\0\\0\\0\\0\\0\\0\\0\\0\\0\\0\\0\\0\\0\\0\\0\\0\\0\\0\\0\\1\\2\\8\\21\\48\\107\\233\\477\\965\\1908\\3726\\7191\\13850\\26553\\51219\\99393\\195289\\389344\\790821\\1636517\\3453698\\7417316\\16176919\\35708318\\79529997\\
\end{array}$
$\begin{array}{r}
\chi_{132}\\
0\\0\\0\\0\\0\\0\\0\\0\\0\\0\\0\\0\\0\\0\\0\\0\\0\\0\\0\\0\\0\\0\\0\\0\\0\\0\\0\\0\\0\\0\\0\\0\\0\\0\\0\\2\\21\\70\\273\\813\\2459\\6796\\18662\\49016\\127573\\323948\\814191\\2014235\\4934162\\11944104\\28648476\\68038799\\
\end{array}$
$\begin{array}{r}
\chi_{133}\\
0\\0\\0\\0\\0\\0\\0\\0\\0\\0\\0\\0\\0\\0\\0\\0\\0\\0\\0\\0\\0\\0\\0\\0\\0\\0\\0\\1\\5\\13\\34\\80\\179\\376\\778\\1550\\3050\\5883\\11264\\21359\\40459\\76539\\145651\\279386\\543062\\1072348\\2158354\\4431005\\9279882\\19789644\\42864534\\93989745\\
\end{array}$
$\begin{array}{r}
\chi_{134}\\
0\\0\\0\\0\\0\\0\\0\\0\\0\\0\\0\\0\\0\\0\\0\\0\\0\\0\\0\\0\\0\\0\\0\\0\\0\\0\\0\\2\\5\\15\\37\\86\\185\\397\\799\\1599\\3111\\6005\\11410\\21651\\40795\\77189\\146504\\281010\\545515\\1077696\\2168327\\4454003\\9329596\\19905755\\43126773\\94597292\\
\end{array}$
$\begin{array}{r}
\chi_{135}\\
0\\0\\0\\0\\0\\0\\0\\0\\0\\0\\0\\0\\0\\0\\0\\0\\0\\0\\0\\0\\0\\0\\0\\0\\0\\0\\0\\0\\0\\0\\0\\0\\0\\0\\0\\0\\2\\23\\108\\446\\1538\\4883\\14552\\41321\\113121\\300750\\780400\\1984152\\4958092\\12201433\\29625584\\71062047\\
\end{array}$
\newpage
$\begin{array}{r}
\chi_{137}\\
0\\0\\0\\0\\0\\0\\0\\0\\0\\0\\0\\0\\0\\0\\0\\0\\0\\0\\0\\0\\0\\0\\0\\0\\0\\0\\0\\0\\0\\0\\2\\6\\18\\52\\127\\303\\697\\1551\\3386\\7313\\15618\\33272\\70932\\151753\\326848\\709895\\1555789\\3440224\\7669212\\17209242\\38807117\\87777806\\
\end{array}$
$\begin{array}{r}
\chi_{138}\\
0\\0\\0\\0\\0\\0\\0\\0\\0\\0\\0\\0\\0\\0\\0\\0\\0\\0\\0\\0\\0\\0\\0\\0\\0\\0\\0\\0\\0\\0\\0\\0\\0\\0\\2\\6\\36\\102\\339\\965\\2779\\7542\\20429\\53562\\139108\\354530\\893243\\2219303\\5454816\\13253442\\31881329\\75927421\\
\end{array}$
$\begin{array}{r}
\chi_{139}\\
0\\0\\0\\0\\0\\0\\0\\0\\0\\0\\0\\0\\0\\0\\0\\0\\0\\0\\0\\0\\0\\0\\0\\0\\0\\0\\0\\0\\0\\4\\9\\27\\62\\155\\332\\738\\1521\\3197\\6471\\13243\\26657\\54346\\110570\\228303\\475324\\1006211\\2156657\\4693584\\10331723\\22989424\\51525740\\116110182\\
\end{array}$
$\begin{array}{r}
\chi_{140}\\
0\\0\\0\\0\\0\\0\\0\\0\\0\\0\\0\\0\\0\\0\\0\\0\\0\\0\\0\\0\\0\\0\\0\\0\\0\\0\\0\\1\\4\\14\\36\\87\\194\\423\\876\\1780\\3528\\6916\\13368\\25735\\49341\\94932\\183653\\359341\\713071\\1440613\\2966514\\6230549\\13327578\\28970516\\63789385\\141808908\\
\end{array}$
$\begin{array}{r}
\chi_{141}\\
0\\0\\0\\0\\0\\0\\0\\0\\0\\0\\0\\0\\0\\0\\0\\0\\0\\0\\0\\0\\0\\0\\0\\0\\0\\0\\0\\0\\0\\0\\6\\12\\43\\101\\253\\553\\1255\\2632\\5640\\11695\\24545\\51038\\107735\\228226\\491599\\1069864\\2363122\\5273174\\11888548\\26972662\\61491807\\140450542\\
\end{array}$
$\begin{array}{r}
\chi_{142}\\
0\\0\\0\\0\\0\\0\\0\\0\\0\\0\\0\\0\\0\\0\\0\\0\\0\\0\\0\\0\\0\\0\\0\\0\\0\\0\\0\\0\\1\\5\\14\\44\\105\\249\\551\\1204\\2523\\5265\\10755\\21949\\44544\\90866\\186158\\385951\\809244\\1722492\\3717460\\8136218\\18013600\\40272034\\90647085\\204960846\\
\end{array}$
$\begin{array}{r}
\chi_{143}\\
0\\0\\0\\0\\0\\0\\0\\0\\0\\0\\0\\0\\0\\0\\0\\0\\0\\0\\0\\0\\0\\0\\0\\0\\0\\0\\0\\0\\3\\9\\26\\67\\161\\354\\775\\1611\\3316\\6673\\13373\\26516\\52818\\105290\\212263\\432579\\896407\\1887395\\4045892\\8807841\\19448096\\43405831\\97667985\\220857196\\
\end{array}$
\newpage
$\begin{array}{r}
\chi_{144}\\
0\\0\\0\\0\\0\\0\\0\\0\\0\\0\\0\\0\\0\\0\\0\\0\\0\\0\\0\\0\\0\\0\\0\\0\\0\\0\\0\\0\\0\\0\\0\\0\\5\\16\\64\\170\\475\\1187\\2992\\7174\\17342\\41096\\97872\\231636\\550107\\1303484\\3092214\\7319517\\17300534\\40751346\\95644345\\223412908\\
\end{array}$
$\begin{array}{r}
\chi_{145}\\
0\\0\\0\\0\\0\\0\\0\\0\\0\\0\\0\\0\\0\\0\\0\\0\\0\\0\\0\\0\\0\\0\\0\\0\\0\\0\\0\\0\\2\\6\\19\\48\\119\\264\\585\\1227\\2560\\5207\\10607\\21411\\43630\\89336\\185845\\392101\\843066\\1842378\\4092747\\9206188\\20920146\\47845726\\109849023\\252481498\\
\end{array}$
$\begin{array}{r}
\chi_{146}\\
0\\0\\0\\0\\0\\0\\0\\0\\0\\0\\0\\0\\0\\0\\0\\0\\0\\0\\0\\0\\0\\0\\0\\0\\0\\0\\0\\0\\0\\1\\2\\11\\28\\81\\191\\473\\1066\\2449\\5414\\12123\\26827\\60297\\135727\\310066\\712752\\1655617\\3863461\\9059551\\21262382\\49902704\\116838212\\272667843\\
\end{array}$
x$\begin{array}{r}
\chi_{147}\\
0\\0\\0\\0\\0\\0\\0\\0\\0\\0\\0\\0\\0\\0\\0\\0\\0\\0\\0\\0\\0\\0\\0\\0\\0\\0\\0\\0\\0\\0\\0\\0\\0\\5\\20\\74\\241\\723\\2031\\5521\\14510\\37365\\94625\\236593\\585363\\1435838\\3494465\\8443650\\20262072\\48292644\\114324172\\268807682\\
\end{array}$
$\begin{array}{r}
\chi_{148}\\
0\\0\\0\\0\\0\\0\\0\\0\\0\\0\\0\\0\\0\\0\\0\\0\\0\\0\\0\\0\\0\\0\\0\\0\\0\\0\\0\\0\\0\\0\\0\\0\\0\\0\\7\\27\\127\\401\\1319\\3824\\11000\\29947\\80583\\210386\\542760\\1374466\\3441750\\8503628\\20797570\\50315286\\120597234\\286334504\\
\end{array}$
$\begin{array}{r}
\chi_{149}\\
0\\0\\0\\0\\0\\0\\0\\0\\0\\0\\0\\0\\0\\0\\0\\0\\0\\0\\0\\0\\0\\0\\0\\0\\0\\0\\0\\0\\0\\0\\0\\5\\15\\53\\142\\392\\953\\2340\\5439\\12740\\29181\\67365\\154631\\357971\\829424\\1934183\\4518169\\10585617\\24795878\\58052983\\135569926\\315598916\\
\end{array}$
$\begin{array}{r}
\chi_{150}\\
0\\0\\0\\0\\0\\0\\0\\0\\0\\0\\0\\0\\0\\0\\0\\0\\0\\0\\0\\0\\0\\0\\0\\0\\0\\0\\0\\0\\0\\0\\3\\9\\32\\79\\211\\491\\1165\\2623\\5971\\13345\\30234\\68435\\156988\\362293\\844661\\1979272\\4662394\\10999332\\25960621\\61153694\\143630088\\335885737\\
\end{array}$
\newpage
$\begin{array}{r}
\chi_{151}\\
0\\0\\0\\0\\0\\0\\0\\0\\0\\0\\0\\0\\0\\0\\0\\0\\0\\0\\0\\0\\0\\0\\0\\0\\0\\0\\0\\0\\0\\0\\5\\15\\50\\131\\338\\793\\1853\\4136\\9207\\20132\\44146\\96531\\212859\\471897\\1057139\\2386816\\5436209\\12455417\\28677871\\66185043\\152880767\\352781947\\
\end{array}$
$\begin{array}{r}
\chi_{152}\\
0\\0\\0\\0\\0\\0\\0\\0\\0\\0\\0\\0\\0\\0\\0\\0\\0\\0\\0\\0\\0\\0\\0\\0\\0\\0\\0\\0\\0\\0\\3\\9\\32\\79\\211\\492\\1176\\2648\\6060\\13562\\30805\\69720\\160032\\368845\\859204\\2010165\\4729027\\11140427\\26263447\\61800036\\145025538\\338901911\\
\end{array}$
$\begin{array}{r}
\chi_{153}\\
0\\0\\0\\0\\0\\0\\0\\0\\0\\0\\0\\0\\0\\0\\0\\0\\0\\0\\0\\0\\0\\0\\0\\0\\0\\0\\0\\0\\0\\0\\0\\0\\3\\12\\50\\150\\458\\1231\\3317\\8448\\21476\\53192\\131452\\320830\\781149\\1888473\\4550849\\10905330\\26019761\\61732274\\145678353\\341691746\\
\end{array}$
$\begin{array}{r}
\chi_{154}\\
0\\0\\0\\0\\0\\0\\0\\0\\0\\0\\0\\0\\0\\0\\0\\0\\0\\0\\0\\0\\0\\0\\0\\0\\0\\0\\0\\0\\0\\0\\0\\0\\4\\16\\63\\184\\544\\1443\\3798\\9539\\23814\\58271\\142097\\343417\\828280\\1988057\\4760382\\11351856\\26974239\\63796113\\150165108\\351529414\\
\end{array}$
$\begin{array}{r}
\chi_{155}\\
0\\0\\0\\0\\0\\0\\0\\0\\0\\0\\0\\0\\0\\0\\0\\0\\0\\0\\0\\0\\0\\0\\0\\0\\0\\0\\0\\0\\0\\1\\3\\11\\31\\79\\195\\464\\1078\\2479\\5698\\13081\\30325\\70686\\166397\\393966\\938588\\2241331\\5360189\\12802182\\30504611\\72393427\\170989044\\401591136\\
\end{array}%
$

\end{document}